\documentclass[11pt]{article}

\usepackage{geometry}                % See geometry.pdf to learn the layout options. There are lots.
\usepackage{graphicx}
\usepackage{amssymb}
\usepackage{amsthm}
\usepackage{amsfonts}
\usepackage{amsmath}
\usepackage{mathrsfs}
\usepackage{amstext}

\usepackage{epstopdf}
\title{Compressed Sensing and Redundant Dictionaries}
\author{Holger Rauhut\footnotemark[1], Karin Schnass and Pierre Vandergheynst\footnotemark[2]}
\newcommand\dico{\mathbf{\Phi}}
\newcommand\atom{\varphi}

\newcommand\ip[2]{\langle #1, #2\rangle}
\newcommand\natoms{K}
\newcommand\sparsity{S}
\newcommand\spacedim{d}
\newcommand\samples{n}

\newcommand\good{\Lambda}
\newcommand\bad{{\overline{\good}}}
\newcommand\sensing{\mathbf{\Psi}}

\newcommand\eps{{\varepsilon}}
\newcommand\nlevel{{\eta}}
\newcommand\ie{{i.e., }}

\newcommand{\N}{{\mathbb{N}}}
\newcommand{\R}{{\mathbb{R}}}

\renewcommand{\P}{{\mathbb{P}}}
\newcommand{\E}{{\mathbb{E}}}

\theoremstyle{plain}
\newtheorem{Theorem}{Theorem}[section]
\newtheorem{lemma}[Theorem]{Lemma}
\newtheorem{corollary}[Theorem]{Corollary}
\newtheorem{example}[Theorem]{Example}

\theoremstyle{remark}

\newenvironment{Proof}{\noindent
{\bf\underline{Proof:} }}
{\hspace*{\fill}\qed\vskip1em}

\graphicspath{figures/}

\numberwithin{equation}{section}

\begin{document}
\maketitle
\begin{abstract}
This article extends the concept of {\it compressed sensing} 
to signals that are
not sparse in an orthonormal basis but rather in a redundant dictionary. 
It is shown that a matrix, 
which is a composition of a random matrix
of certain type and a deterministic dictionary,
has small restricted isometry constants. Thus, signals that are
sparse with respect to the dictionary can be recovered via Basis Pursuit from
a small number of random measurements. 
Further, thresholding is investigated as 
recovery algorithm for compressed sensing and
conditions are provided that guarantee reconstruction with high probability.
The different schemes are compared by numerical experiments.
\end{abstract}


%\noindent
%{\bf AMS  Subject Classification}

\renewcommand{\thefootnote}{\fnsymbol{footnote}}

\footnotetext[1]{Numerical Harmonic Analysis Group, 
Faculty of Mathematics, University of Vienna,\\
Nordbergstrasse 15, A-1090 Vienna, Austria,
%email: 
{\tt holger.rauhut@univie.ac.at}\\ 
HR acknowledges the financial support provided by the
European Union's Human Potential Programme under contract MEIF-CT-2006-022811.}
\footnotetext[2]{Signal Processing Institute, Swiss Federal Institute of Technology (EPFL), CH-1015 Lausanne, Switzerland, {\tt name.surname@epfl.ch}}
\renewcommand{\thefootnote}{\arabic{footnote}}

\section{Introduction}

Recently there has been a growing interest in recovering sparse 
signals from their projection 
onto a small number of random vectors
\cite{carota06-1,cataXX,do04,gitr05,ra05-7,ru06-1}. 
The word most often used 
in this context 
is {\it compressed sensing}. It originates from the idea that
it is not necessary to invest a lot of power into 
observing the entries of a sparse signal
in all coordinates when most of them are zero anyway. Rather it should be 
possible 
to collect only a small number of measurements that still allow for reconstruction. This is potentially useful in applications where
one cannot afford to collect or transmit a lot of measurements but has rich resources at the decoder.

Until now the theory of compressed sensing has only been developed for 
classes of 
signals that have a very sparse representation in an orthonormal basis (ONB). 
This is a rather stringent restriction. Indeed, allowing the signal
to be sparse with respect to a redundant dictionary adds a lot of flexibility
and significantly extends the range of applicability.
Already the use of two ONBs instead of 
just one dramatically increases the class of signals that can be 
modelled in this way. 
A more practical example would be a dictionary made up of damped sinusoids which is used for NMR spectroscopy, see \cite{dr06}.

Before we can go into further explanations about the scope of this paper it is necessary to provide some
background information.
The basic problem in compressed sensing is to determine the minimal
number $\samples$ of linear non-adaptive measurements that allows for (stable)
reconstruction of a signal $x \in \R^\spacedim$ that has at most
$\sparsity$ non-zero components.
Additionally, one
requires that this task can be performed reasonably fast.
Each of the $\samples$ measurements can be written as an inner product 
of the sparse signal $x \in \R^\spacedim$ 
with a vector in $\R^\spacedim$. 
To simplify the notation we store all the 
vectors as rows in a matrix $\sensing \in \R^{\samples \times \spacedim}$
and all the measurements in the $\samples$-dimensional
vector $s=\sensing x$. 

A naive approach to the problem of recovering $x$ from $s$ 
consists in solving
the $\ell_0$ minimization problem 
\begin{equation}\nonumber
(P_0) \hspace{1cm} \min \|x\|_0 \mbox{ subject to } \|s-\sensing x\|_2\leq\nlevel,
\end{equation}
where $\nlevel$ is the expected noise on the measurements, 
$\|\cdot\|_0$ counts the number of  non-zero entries of $x$ and 
$\|\cdot\|_2$ denotes the standard Euclidean norm.
Although there are simple recovery conditions available, the above
approach is not reasonable in practice because its solution
is NP-hard \cite{avdama97,na95}.

In order to avoid this severe drawback there have been basically two 
approaches proposed in the signal recovery community. 
The first is using greedy algorithms like {\it Thresholding} \cite{grmarascva06}
or {\it (Orthogonal) Matching Pursuit} (OMP) \cite{mazh93,tr04}. 
Thresholding simply 
calculates the inner products of the signal with all atoms, finds 
the ones with largest absolute values and then calculates the 
orthogonal projection onto the span of the corresponding atoms. 
OMP works iteratively by picking the atoms in a greedy fashion. 
In each step it finds the atom with highest absolute inner product with the 
residual and adds it to the already found atoms. Then it calculates a 
new approximant by projecting the signal on the linear span of 
the already found atoms and a new residual by subtracting the 
approximant from the signal, cp.~Table~\ref{table:greedy}.

\begin{table}[bdp]\label{table:greedy}
\caption{Greedy Algorithms}
Goal: reconstruct $x$ from $s = \sensing x$\\ 
columns of $\sensing$ denoted by $\psi_j$, $\sensing_\good^\dagger$: pseudo-inverse of $\sensing_\good$\\
\begin{tabular}{|l|l|}
\hline
{\bf OMP} & {\bf Thresholding} \\ \hline
initialise:  $z=0$, $r=s$, $\good=\emptyset$ & find: $\good$ that contains the indices\\
find:   $i=\arg \max_j |\ip{r}{\psi_j}|$     & \phantom{find: }corresponding to the $\sparsity$ largest\\
update: $\good=\good \cup \{i\}$, $r=s-\sensing_\good \sensing_\good^\dagger s$ & 
\phantom{find: }values of $|\ip{s}{\psi_j}|$\\
iterate until stopping criterion is attained & output: $x=\sensing_\good^\dagger s$\\
output: $x = \sensing^\dagger_\good s$ & \\ \hline
\end{tabular}
\end{table}

%%%%%%% for two paragraphs

%\begin{table}[htbdp]\label{table:greedy}
%\caption{Greedy Algorithms}
%Goal: reconstruct $x$ from $s = \sensing x$\\ 
%columns of $\sensing$ denoted by $\psi_j$,\\
%$\sensing_\good^\dagger$: pseudo-inverse of $\sensing_\good$\\
%\begin{tabular}{|l|}
%\hline
%{\bf OMP} \\ \hline
%initialise:  $z=0$, $r=s$, $\good=\emptyset$ \\
%find:   $i=\arg \max_j |\ip{r}{\psi_j}|$    \\
%update: $\good=\good \cup \{i\}$, $r=s-\sensing_\good \sensing_\good^\dagger s$ \\
%iterate until stopping criterion is attained \\
%output: $x = \sensing^\dagger_\good s$ \\ \hline
% {\bf Thresholding} \\ \hline
%find: $\good$ that contains the indices\\
% \phantom{find: }corresponding to the $\sparsity$ largest\\ 
%\phantom{find: }values of $|\ip{s}{\psi_j}|$\\ 
% output: $x=\sensing_\good^\dagger s$\\ \hline
%\end{tabular}
%\end{table}

The second approach is the {\it Basis Pursuit} (BP) principle.
Instead of considering $(P_0)$ one solves its convex relaxation
\begin{equation}\nonumber%\label{BPdn}
(P_1) \hspace{1cm}\min \|x\|_1 \mbox{ subject to } \|s-\sensing x\|_2<\nlevel,
\end{equation}
where $\|x\|_1 = \sum |x_i|$ denotes the $\ell_1$-norm.
This can be done via linear programming in the real case and 
via cone programming in the complex case.
Clearly, one hopes that the solutions of $(P_0)$ and $(P_1)$ coincide, 
see \cite{chdosa99, doelte06} for details.

Both approaches pose certain requirements on the 
matrix $\sensing$ in order to ensure recovery success.
Recently, Cand{\`e}s, Romberg and Tao \cite{carota06-1,cataXX} observed
that successful recovery by BP is guaranteed whenever
$\sensing$ obeys a uniform uncertainty principle.
Essentially this means that every submatrix of $\sensing$ of a certain
size has to be well-conditioned. More precisely, 
let $\good \subset \{1,\hdots,\spacedim\}$ and $\sensing_\good$ be the submatrix
of $\sensing$ consisting of the columns indexed by $\good$. The local
isometry constant $\delta_\good = \delta_\good(\sensing)$ is the smallest number
satisfying
\begin{equation}\label{def:local_iso}
(1-\delta_\good) \|x\|_2^2 \leq \|\sensing_\good x\|_2^2 \leq (1+\delta_\good)\|x\|_2^2,
\end{equation}
for all coefficient vectors $x$ supported on $\good$. The (global) restricted 
isometry constant is then defined as
\[
\delta_\sparsity = \delta_\sparsity(\sensing) \,:=\, \sup_{|\good| = \sparsity}
\delta_\good(\sensing),\quad \sparsity \in \N.
\]
The matrix $\sensing$ is said to satisfy a uniform uncertainty principle
if it has small restricted isometry constants, say 
$\delta_\sparsity(\sensing) \leq 1/2$. Based on this concept, Cand{\`e}s, Romberg and Tao proved the following
recovery theorem for BP in \cite[Theorem 1]{carota06-1}.

\begin{Theorem}\label{thm:BP}
Assume that $\sensing$ satisfies
\[
\delta_{3\sparsity}(\sensing) + 3 \delta_{4\sparsity}(\sensing) < 2
\]
for some $\sparsity \in \N$. Let $x$ be an $\sparsity$-sparse vector 
and assume we are given noisy 
data $y=\sensing x+\xi$ with $\|\xi\|_2 \leq \nlevel$. 
Then the solution $x^\#$ 
to the problem $(P_1)$
satisfies
\begin{equation}\label{BP_error}
\|x^\# - x\|_2 \,\leq\, C \nlevel. % + C_2 \frac{\|x-x_M\|_1}{\sqrt{M}}.
\end{equation}
The constant $C$ depends only on $\delta_{3\sparsity}$ and 
$\delta_{4\sparsity}$. 
If $\delta_{4\sparsity} \leq 1/3$ then $C \leq 15.41$. 
\end{Theorem}
In particular, if no noise is present, i.e., $\nlevel = 0$, then 
under the stated condition BP recovers $x$ exactly.
Note that a slight variation of the above theorem
holds also in the case that $x$ is not sparse in a strict sense,
but can be well-approximated by an $\sparsity$-sparse vector 
\cite[Theorem 2]{carota06-1}. 

Of course, the above theorem is only useful if there are matrices
$\sensing$ satisfying the uniform uncertainty principle. So far 
no deterministic
construction is known (for a reasonably small ratio $\samples / \sparsity$).
However, an $\samples \times \spacedim$ random matrix 
with entries drawn from a standard Gaussian distribution 
(or some other distribution showing certain 
concentration properties, see below) 
will have small restricted isometry constants $\delta_\sparsity$ 
with 'overwhelming probability' as long as
\begin{equation}
\samples = {\cal O}(\sparsity \log(\spacedim/\sparsity)),
\end{equation}
see \cite{badadewa06, carota06-1, cataXX, ru06-1} for details.

The results for OMP in compressed sensing are weaker than for BP.  
While it can again be shown that with high probability a signal 
can be reconstructed from the random measurements 
$\sensing x$ if $\samples > C \sparsity \log{\spacedim}$, this result is 
no longer uniform in the sense that no single measurement matrix $\sensing$ will 
simultaneously work for all possible sparse signals, see \cite{gitr05}.
 
As already announced we want to address the question
whether the techniques described above can 
be extended to signals $y$ that are not sparse in an ONB but rather in a 
redundant dictionary $\dico \in \R^{\spacedim \times \natoms}$ with 
$\natoms > \spacedim$. 
So now $y=\dico x$, where $x$ has only few non-zero 
components. Again the goal is to reconstruct $y$ from few measurements.
More formally, given a suitable measurement matrix 
$A \in \R^{\samples \times \spacedim}$ we want to recover
$y$ from $s = Ay = A \dico x$. 
The key idea then is to use the 
sparse representation 
in $\dico$ to drive the reconstruction procedure, 
\ie try to identify the sparse coefficient 
sequence $x$ and from that reconstruct $y$.
Clearly, we may represent $s = \sensing x$ with 
\[
\sensing = A \dico \in 
\R^{\samples \times \natoms}.
\]
In particular, we can apply all of the reconstruction methods
described above by using this particular matrix $\sensing$.
Of course, the remaining question is whether for a fixed %(deterministic)
dictionary $\dico \in \R^{\spacedim \times \natoms}$ one can find a suitable matrix 
$A \in \R^{\samples \times \spacedim}$ such that the composed matrix 
$\sensing = A \dico$ allows for reconstruction of vectors having 
only a small number of non-zero
entries. Again the strategy is to choose a random 
matrix $A$, for instance with independent
standard Gaussian entries, and investigate under which conditions 
on $\dico$, $\samples$ and $\sparsity$ recovery is successful with
high probability.

Note that already Donoho considered extensions
from orthonormal bases to (redundant) tight frames $\dico$ in \cite{do04}. There it is assumed that the analysis 
coefficients $x' = \dico^* y = \dico^* \dico x$
are sparse. For redundant frames, however, this assumption does not seem very realistic as even for
sparse vectors $x$ the coefficient vector $x' = \dico^* \dico x$ is usually fully populated.

In the following section we will investigate under which conditions on the deterministic 
dictionary $\dico$ its combination with a random measurement matrix will have small isometry 
constants. By Theorem \ref{thm:BP} this determines 
how many measurements $\samples$ will be typically required for BP to succeed in reconstructing 
all signals of sparsity $\sparsity$ with respect to the given dictionary. 
In Section \ref{Sec:Thresh} we will analyse the performance of thresholding, 
which actually has not yet been considered as a 
reconstruction algorithm in compressed sensing
because of its simplicity and hence resulting limitations.
The last section is dedicated to numerical simulations showing the performance of compressed 
sensing for dictionaries in practice and comparing it to the situation where sparsity is 
induced by an ONB. Even though we have not yet been able to theoretically analyse OMP 
for compressed sensing we will do simulations for all three approaches.

%%% Local Variables:
%%% TeX-master: "CompSensing.tex"
%%% End:

\section{Isometry Constants for $A\dico$}

In order to determine the isometry constants for a matrix of the 
type $\sensing = A\dico$, where $A$ is an $\samples \times \spacedim$ 
measurement matrix
and $\dico$ is a $\spacedim \times \natoms$ dictionary, we will 
follow the approach taken in \cite{badadewa06}, which was inspired by proofs 
for the Johnson-Lindenstrauss lemma \cite{ac01}. We will not discuss this 
connection further but use as starting point concentration of measure 
for random variables. This describes the phenomenon 
that in high dimensions the probability mass of certain random variables 
concentrates strongly around their expectation. 

In the following we will assume that $A$ is an $\samples \times \spacedim$
random matrix that satisfies 
\begin{equation} \label{eq:conc}
\P\left(\big| \|Av\|^2 - \|v\|^2 \big| \geq \eps \|v\|^2\right) 
\leq 2 e^{-c \frac{n}{2}\eps^2}, \hspace{2em} \eps \in (0,1/3)
\end{equation}
for all $v \in \R^\spacedim$ and some constant $c > 0$.
Let us list some examples of random matrices that satisfy the above
condition.
\begin{itemize}
\item {\bf Gaussian ensemble:} If the entries of $A$ are independent
normal variables with mean zero and variance $\samples^{-1}$ then
\begin{equation} \label{eq:Gauss_conc}
\P(\big| \|Av\|^2 - \|v\|^2 \big| \geq \eps \|v\|^2) \leq 2
e^{-\frac{n}{2}(\eps^2/2 -\eps^3/3)}, \hspace{2em} \eps \in (0,1),
\end{equation}
see e.g.~\cite{ac01,badadewa06}.
In particular, (\ref{eq:conc}) holds with $c = 1/2 - 1/9 = 7/18$. 
%for
%$\eps_0 \leq 1/3$. 
\item {\bf Bernoulli ensemble:} Choose the entries of $A$ as independent
realizations of $\pm 1/\sqrt{\samples}$ random variables. Then again 
(\ref{eq:Gauss_conc}) is valid, see \cite{ac01,badadewa06}. In 
particular (\ref{eq:conc}) holds with $c = 7/18$. % for $\eps_0 \leq 1/3$.
\item {\bf Isotropic subgaussian ensembles:} In generalization of the
two examples above, we can choose the rows of $A$ as 
$\frac{1}{\sqrt{\samples}}$-scaled independent copies of a random vector 
$Y\in \R^\spacedim$ that satisfies $\E |\langle Y,v\rangle|^2 = \|v\|^2$ for all
$v \in \R^\spacedim$
and has subgaussian tail behaviour. See \cite[eq.~(3.2)]{mepato06} for details.
\item {\bf Basis transformation:} If we take any valid random matrix $A$ 
and a (deterministic) orthogonal $\spacedim \times \spacedim$ 
matrix $U$ then it is easy to see that also $AU$ satisfies the concentration
inequality (\ref{eq:conc}). In particular, this applies to the Bernoulli
ensemble although in general $AU$ and $A$ have different probability distributions.
\end{itemize}

Using the concentration inequality (\ref{eq:conc}) 
we can now investigate the 
local and subsequently the global restricted isometry constants of  
the $n \times \natoms$ matrix $A\dico$.

\begin{lemma}\label{lemma:locisometry}
Let $A$ be a random matrix of size $\samples \times \spacedim$
drawn from a distribution that satisfies the concentration
inequality~\eqref{eq:conc}. Extract from the 
$\spacedim \times \natoms$ dictionary $\dico$ 
any sub-dictionary
$\dico_\good$ of size $\sparsity$, \ie $|\good|=\sparsity$ with 
(local) isometry constant $\delta_\good = \delta_\good(\dico)$. For $0<\delta<1$
we set 
\begin{equation}\label{eq:nu}
\nu := \delta_\good + \delta +\delta_\good \delta. 
\end{equation}
Then
\begin{equation}\label{localisometry}
(1 - \nu)\|x\|^2\leq \|A\dico_\good x\|^2 \leq \|x\|^2 (1 + \nu)
\end{equation}
with probability exceeding
\begin{equation}\label{Plocisometry}
1- 2 \left(1+\frac{12}{\delta}\right)^\sparsity 
e^{- \frac{c}{9} \delta^2 \samples}.
\end{equation}
\end{lemma}

\begin{Proof}
First we choose a finite $\eps_1$-covering of the unit sphere in
$\R^\sparsity$, \ie a set of points $Q$, with $\|q\|= 1$ for all
$q \in Q$, such that for all $\|x\|=1$
\begin{eqnarray*}
\min_{q \in Q}\|x-q \|\leq \eps_1 %,\qquad 0<\eps_1 < 1.
\end{eqnarray*}
for some $\eps_1 \in (0,1)$.
According to Lemma~2.2 in \cite{mepato06} there exists such a $Q$ with
$|Q|\leq (1 + 2/\eps_1)^\sparsity$. Applying the measure
concentration in \eqref{eq:conc} with $\eps_2 < 1/3$ to all the points $\dico_\good q$ and taking the union bound we get
\begin{equation}
(1 - \eps_2)\|\dico_\good q\|^2\leq \|A\dico_\good q\|^2 
\leq (1 +\eps_2) \|\dico_\good q\|^2 \hspace{1em} \mbox{ for all } 
q \in Q,
\end{equation}
with probability larger than 
\[
1- 2 \left(1+\frac{2}{\eps_1}\right)^\sparsity e^{-c n \eps_2^2}.
\]
Define $\nu$ as the smallest number such that
\begin{eqnarray}\label{definemu}
\|A\dico_\good x\|^2\leq (1 + \nu) \|x\|^2, \hspace{1em} \mbox{
for all } x \mbox{ supported on }\good.
\end{eqnarray}
Now we estimate $\nu$ in terms of $\eps_1,\eps_2$. We know that for all
$x$ with $\|x\|=1$ we can choose a $q$ such that $\|x-q\|\leq
\eps_1$ and get
\begin{eqnarray*}
\|A\dico_\good x \| &\leq&  \|A\dico_\good q\|+\|A\dico_\good
(x-q)\|\\
&\leq&  (1+\eps_2)^{\frac{1}{2}} \|\dico_\good q\|+\|A\dico_\good
(x-q)\|\\
&\leq&  (1+\eps_2)^{\frac{1}{2}} (1+\delta_\good)^{\frac{1}{2}}
+(1+\nu)^{\frac{1}{2}}\eps_1.
\end{eqnarray*}
Since $\nu$ is the smallest possible constant for which
(\ref{definemu}) holds it also has to satisfy
\begin{eqnarray*}
\sqrt{1+\nu} &\leq&  \sqrt{1+\eps_2}
\sqrt{1+\delta_\good} +\eps_1\sqrt{1+\nu}.
\end{eqnarray*}
Simplifying the above equation yields
\[
(1+\nu) \leq \frac{1+\eps_2}{(1-\eps_1)^2}(1+\delta_\good).
\]
Now we choose $\eps_1=\delta/6$ and $\eps_2=\delta/3<1/3$.
Then
\begin{align}
\frac{1+\eps_2}{(1-\eps_1)^2} = & \frac{1+\delta/3}{(1-\delta/ 6)^2}=\frac{1+\delta/3}{1-\delta/3 + \delta^2/36}
< \frac{1+\delta/3}{1-\delta/3} = 1 + \frac{2\delta/3}{1-\delta/3} < 
1 +  \delta.
\notag
\end{align} 
Thus,
\begin{eqnarray*}
\nu& < & \delta + \delta_\good(1 + \delta).
\end{eqnarray*}
To get the lower bound we operate in a similar fashion.
\begin{eqnarray*}
\|A\dico_\good x \| &\geq&  \|A\dico_\good q\|- \|A\dico_\good
(x-q)\| \\
&\geq&  (1-\eps_2)^{\frac{1}{2}} (1-\delta_\good)^{\frac{1}{2}}
-(1+\nu)^{\frac{1}{2}}\eps_1.
\end{eqnarray*}
Now square both sides and observe that $\nu<1$ (otherwise we have
nothing to show). Then we finally arrive at
\begin{eqnarray*}
\|A\dico_\good x \|^2 &\geq& \big((1-\eps_2)^{\frac{1}{2}}
(1-\delta_\good)^{1/2} -\eps_1\sqrt{2}\big)^2  \\&\geq&
(1-\eps_2) (1-\delta_\good) -2\eps_1\sqrt{2}\sqrt{1-\eps_2}\sqrt{1-\delta_\good}
+ 2\eps_1^2\\&\geq& 1 - \delta_\good - \eps_2 - 2 \eps_1 \sqrt{2}
\geq 1 - \delta_\good - \delta \geq 1 - \nu.
\end{eqnarray*}
This completes the proof.
\end{Proof}
Note that the choice of $\eps_1$ and $\eps_2$ in the previous proof 
is not the only one possible. While our choice has the advantage of resulting in an appealing form of $\nu$ in (\ref{eq:nu}), others might actually yield better constants.

Based on the previous theorem it is easy to derive an estimation
of the global restricted 
isometry constants of the composed matrix $\sensing = A \dico$.

\begin{Theorem}\label{thm:uniformisometry}
Let $\dico \in \R^{\spacedim \times \natoms}$ be a dictionary 
of size $\natoms$ in $\R^\spacedim$ with 
restricted isometry constant
$\delta_\sparsity(\dico)$, $S \in \N$. Let $A \in \R^{\samples \times \spacedim}$
be a random matrix satisfying (\ref{eq:conc}) and assume
\begin{equation}\label{ineq:samples}
\samples \geq C \delta^{-2}\left(\sparsity \log(K/\sparsity) 
+ \log(2e(1+12/\delta)) + t \right)
\end{equation}
for some $\delta \in (0,1)$ and $t > 0$. Then with probability at 
least $1-e^{-t}$ 
the composed matrix $\sensing = A\dico$ has restricted isometry constant
\begin{equation}\label{delta:composed}
\delta_\sparsity(A\dico)\leq \delta_\sparsity(\dico) + \delta(1+ \delta_\sparsity(\dico)).
\end{equation}
The constant satisfies $C \leq 9/c$.
\end{Theorem}
\begin{Proof}
By Lemma~\ref{lemma:locisometry} we can estimate
the probability that a sub-dictionary $\sensing_\good = (A\dico)_\good=A\dico_\good$, 
$\good \subset \{1,\hdots,K\}$ fails 
to have (local) isometry constants
$
\delta_\good(\sensing) \leq \delta_\good(\dico) + \delta +\delta_\good(\dico)\delta
$
by
\begin{equation}\nonumber
\P\big(\delta_\good(\sensing) > \delta_\good(\dico) + \delta +\delta_\good(\dico)\delta\big) \leq 2\big(1+\frac{12}{\delta}\big)^\sparsity e^{-\frac{c}{9} \delta^2 n}.
\end{equation} 
By taking the union bound over all ${\natoms \choose \sparsity}$ 
possible sub-dictionaries of size $\sparsity$ we can estimate
the probability of 
$\delta_S(\sensing) = 
\sup_{\good \subset \{1,\hdots,\natoms\},|\good| = \sparsity} \delta_\good(\sensing)$
{\em not} satisfying (\ref{delta:composed}) by
\[
\P\big(\delta_\sparsity(\sensing )> \delta_\sparsity(\dico) + \delta(1+ \delta_\sparsity(\dico))\big) \leq 2 {\natoms \choose \sparsity} \left(1+\frac{12}{\delta}\right)^\sparsity
 e^{-\frac{c}{9} \delta^2 n}.
\]
Using ${\natoms \choose \sparsity} \leq (e \natoms / \sparsity)^{\sparsity}$
(Stirling's formula) and requiring that the above term is less than
$e^{-t}$ shows the claim.
\end{Proof}

Note that for fixed $\delta$ and $t$ 
condition (\ref{ineq:samples}) can be expressed in the more compact
form 
\[
\samples \geq C \sparsity \log(\natoms/ \sparsity).
\]
Moreover, if the dictionary $\dico$ is an orthonormal basis 
then $\delta(\dico) = 0$ and we recover essentially the previously
known estimates of the isometry constants for a random matrix
$A$, see e.g.~\cite[Theorem 5.2]{badadewa06}.  

Now that we have established how the isometry constants of a 
deterministic dictionary $\dico$ are affected by multiplication with a random 
measurement matrix, we only need some more initial information about $\dico$, 
before we can finally apply the result to compressed sensing of signals that 
are sparse in $\dico$.
The following little lemma gives a very crude estimate of the isometry 
constants of $\dico$ in terms of its coherence $\mu$ or Babel 
function $\mu_1(k)$, which are defined as
\begin{equation}\label{def:Babel}
\mu:=\max_{i\neq j} |\ip{\atom_i}{\atom_j}|, \hspace{1cm} \mu_1(k):=\max_{|\good|=k, j \notin \good} \sum_{i\in\good} |\ip{\atom_i}{\atom_j}|.
\end{equation}
\begin{lemma}\label{lem:coh} For a dictionary with coherence $\mu$ 
and Babel function $\mu_1(k)$ we can bound the restricted isometry constants by
\begin{equation}\label{coherence:estimate}
\delta_\sparsity\leq \mu_1(\sparsity -1)\leq (\sparsity -1)\mu.
\end{equation}
\end{lemma}
\begin{Proof} Essentially this can be 
derived from the proof of Lemma~2.3 in \cite{tr04}.
\end{Proof}

Combining this Lemma with Theorem~\ref{thm:uniformisometry} provides
the following estimate of the isometry constants 
of the composed matrix $\sensing = A \dico$.
\begin{corollary}\label{cor:isometry} Let $\dico \in \R^{\spacedim \times \natoms}$ be a dictionary
with coherence $\mu$.
Assume that 
\begin{equation}\label{cond:sparsity}
\sparsity -1 \leq \frac{1}{16}\mu^{-1}.
\end{equation}
Let 
$A \in \R^{\samples\times \spacedim}$ be a random matrix satisfying (\ref{eq:conc}).
Assume that
\[
n \geq C_1 (S \log(K /S) + C_2 + t).
\]
Then with probability at least $1-e^{-t}$ the composed matrix $A \dico$
has restricted isometry constant 
\begin{equation}\label{delta4}
\delta_\sparsity(\sensing) \leq 1/3.
\end{equation}
The constants
satisfy $C_1 \leq 138.51\, c^{-1}$ and 
$C_2 \leq \log(1250/13) + 1 \approx 5.57$. In particular, for the Gaussian and
Bernoulli ensemble $C_1 \leq 356.18$.
\end{corollary}
\begin{Proof} By Lemma \ref{lem:coh} the restricted 
isometry constant of $\dico$ satisfies
\[
\delta_\sparsity(\dico) \leq (S-1)\mu \leq 1/16.
\]
Hence, choosing $\delta = 13/(3\cdot 17)$ yields
\[
\delta(A \dico) \leq \delta_{\sparsity}(\dico) + \delta(1+\delta_{\sparsity}(\dico)) \leq \frac{1}{16} + \frac{13}{3\cdot 17}(1+\frac{1}{16}) = 1/3.
\]
Plugging this particular choice of $\delta$ into Theorem~\ref{thm:uniformisometry} yields the assertion.
\end{Proof}
Of course, the numbers $1/16$ and $1/3$ in (\ref{cond:sparsity}) 
and (\ref{delta4}) were just arbitrarily chosen.
Other choices will only result in different constants $C_1,C_2$.
Combining the previous result with Theorem~\ref{thm:BP} 
yields a result on stable recovery
by Basis Pursuit of sparse signals in a redundant dictionary. 
We leave the straightforward task of formulating
the precise statement to the interested reader. 
We just want to point out that this recovery result is uniform in the sense that a 
single matrix $A$ can ensure recovery of {\em all} 
sparse signals.

The constants $C_1$ and $C_2$ of Corollary \ref{cor:isometry} are
probably not optimal. In the case of a Gaussian ensemble $A$ and
an orthonormal basis $\dico$ recovery conditions for BP with quite
small constants were obtained in \cite{ru06-1} and precise asymptotic results can be found in \cite{dota06}.
One might raise the objection that the condition
$S-1 \leq \frac{1}{16\mu}$ in Corollary \ref{cor:isometry} is 
too weak for practial applications. 
A lower bound on the coherence in terms of the dictionary size 
is $$\mu>\sqrt{\frac{\natoms-\spacedim}{\spacedim(\natoms-1)}} $$ and 
for reasonable dictionaries we can usually expect the coherence to 
be of the order $\mu \sim C/\sqrt{d}$. The restriction on the sparsity 
thus is $\sparsity < \sqrt{\spacedim}/C$. 
However, compressed sensing is only useful if indeed the sparsity
is rather small compared to the dimension $\spacedim$, so this restriction
is actually not severe. Moreover, 
if it is already impossible to recover the support 
from complete information on the original signal 
we cannot to expect to do this with even less information.

To illustrate the theorem let us 
have a look at an example where the dictionary is the union of two ONBs.

\begin{example}[Dirac-DCT] Assume that our dictionary is the union of the 
Dirac and the Discrete Cosine Transform bases in $\R^d$ for $d=2^{2p+1}$. 
The coherence in this case is $\mu=\sqrt{2/d}=2^{-p}$ and the number of 
atoms $K= 2^{2p+2}$. If we assume the sparsity of the signal to be smaller 
than $2^{p-6}$ we 
get the following crude estimate for the number of necessary samples to 
have $\delta_{4\sparsity}(A\dico)<1/3$ as recommended for recovery by BP in 
Theorem~\ref{thm:BP},
\begin{equation}\nonumber
\samples \geq  C_1 (4\sparsity(2p \log{2} - \log{\sparsity}) + C_2 + t)
\end{equation}
with the constants $C_1\approx 138.51\,c^{-1}$ and $C_2 \approx 5.57$ 
from Corollary \ref{cor:isometry}.

In comparison if the signal is sparse in just the Dirac basis we can 
estimate the necessary number of samples to have $\delta_{4\sparsity}(A)<1/3$ 
with Theorem~\ref{thm:uniformisometry} as
\begin{equation}\nonumber
\samples \geq C_1' (4\sparsity(2p\log{2}-\log{2\sparsity}) + C_2' + t)
\end{equation}
with $C_1' = \big(\frac{13}{17}\big)^2 C_1$ and $C_2' \approx 5.3$, thus implying
an improvement of roughly the factor $(\frac{17}{13})^2 \approx 1.71$.
\end{example}

%%% Local Variables:
%%% TeX-master: "CompSensing.tex"
%%% End:

\section{Recovery by Thresholding}
\label{Sec:Thresh}

In this section we investigate recovery from random
measurements by thresholding. Since thresholding works by comparing 
inner products of the signal with the atoms
%inside to those of the signal and the atoms outside the support 
an essential ingredient will be stability of inner 
products under multiplication with a random matrix $A$, \ie
\begin{equation}\nonumber
\ip{Ax}{Ay}\approx \ip{x}{y}.
\end{equation}
The exact result that we will use is summarised in the following lemma.

\begin{lemma}\label{lemma:ipconc} Let $x,y \in \R^\spacedim$ with 
$\|x\|_2,\|y\|_2 \leq 1$.
Assume that $A$ is an $\samples \times \spacedim$ random
matrix with independent ${\cal N}(0,\samples^{-1})$ 
entries (independent of $x,y$).
Then for all $t > 0$
%\marginpar{add version without factor 2 for upper and lower seperated}
%\[
%\P\left(\langle A x, Ay\rangle \geq \langle x,y\rangle + t\right)
%\leq \exp\left(-\samples \frac{t^2}{C_1 + C_2 t}\right)
%\]
%and
\begin{align}\label{sec:ineq}
\P\big(|\langle A x, Ay\rangle - & \langle x,y\rangle| \geq  t \big) \notag\\
&\leq 2\exp\left(-\samples \frac{t^2}{C_1 + C_2 t}\right), 
\end{align}
with $C_1 =  \frac{4e}{\sqrt{6\pi}}
\approx 2.5044$ and $C_2 =  e \sqrt{2} \approx 3.8442$. 

The analogue statement holds for a random matrix $A$ with independent
$\pm 1/\sqrt{\samples}$ Bernoulli entries.
\end{lemma}

%The same result (even with the same constants) holds
%also for random matrices $A$ with independent $\pm \samples^{-1}$ entries
%(Bernoulli) matrices.\\

Note that taking $x=y$ in the lemma provides
the concentration inequality (\ref{eq:conc}) for Gaussian
and Bernoulli matrices (with non-optimal constants however).

The proof of the lemma is rather technical and therefore safely 
locked away in Appendix~\ref{app:proof} awaiting inspection by the 
genuinely interested reader there. However armed with it, 
we can now investigate 
the stability of recovery via thresholding.

\begin{Theorem}\label{th:recth} Let $\dico$ be a $\spacedim \times \natoms$
dictionary.
Assume that the support $x$ of a signal $y=\dico x$, normalised to 
have $\|y\|_2=1$, could be recovered by thresholding with a margin $\eps$, \ie
\begin{equation}
\min_{i\in\good}|\ip{y}{\atom_i}|> \max_{k\in\bad}|\ip{y}{\atom_k} | +\eps.\nonumber
\end{equation}
Let $A$ be an $\samples \times \spacedim$ random matrix satisfying one
of the two probability models of the previous lemma.
Then with probability exceeding $1-e^{-t}$ the support and thus 
the signal can be reconstructed via thresholding from the 
$\samples$-dimensional 
measurement vector $s = A y = A\dico x$ as long as
\begin{equation}
\samples \geq C(\eps) (\log{(2\natoms) + t)}. \nonumber
\end{equation}
where $C(\eps) = 4C_1 \eps^{-2} + 2C_2 \eps^{-1}$ and
$C_1, C_2$ are the constants from Lemma \ref{lemma:ipconc}.
In particular, 
\[
C(\eps) \leq C_3 \eps^{-2}
\]
with $C_3 \leq 4C_1 + 2C_2 \leq 17.71$. 
%where $C_3(\eps)<4\eps^{-2}(C_1 + C_2\frac{\eps}{2})<51\eps^{-2}$.
\end{Theorem}

\begin{Proof} Thresholding will succeed if we have
\begin{equation}
\min_{i\in\good}|\ip{Ay}{A\atom_i}|> \max_{k\in\bad}|\ip{Ay}{A\atom_k} |.\nonumber
\end{equation}
So let us estimate the probability that the above inequality does {\em not} 
hold,
\begin{align}
& \P(\min_{i\in\good}|\ip{Ay}{A\atom_i}| \leq \max_{k\in\bad}|\ip{Ay}{A\atom_k} |)\nonumber \\
&\leq  \P(\min_{i\in\good}|\ip{Ay}{A\atom_i}| \leq \min_{i\in\good}|\ip{y}{\atom_i}| -\frac{\eps}{2}) \notag\\
&\hspace{2em}+\P(\max_{k\in\bad}|\ip{Ay}{A\atom_k} | \geq \max_{k\in\bad}|\ip{y}{\atom_k} | +\frac{\eps}{2})\nonumber
\end{align}
The probability of the good components having responses lower than the threshold can be further estimated as
\begin{align}
 \P(\min_{i\in\good}&|\ip{Ay}{A\atom_i}| \leq \min_{i\in\good}|\ip{y}{\atom_i}| -\frac{\eps}{2}) \notag\\
%= 1 - \P(\min_{i\in\good}|\ip{Ay}{A\atom_i}| \leq \min_{i\in\good}|\ip{y}{\atom_i}| -\frac{\eps}{2}) \nonumber \\
&\leq  \P\left(\bigcup_{i\in \good}\{ |\ip{Ay}{A\atom_i}| \leq |\ip{y}{\atom_i}| -\frac{\eps}{2}\}\right)\nonumber \\
&\leq \sum_{i\in\good} \P\left( |\ip{y}{\atom_i} - \ip{Ay}{A\atom_i}| \geq \frac{\eps}{2}\right)\notag\\
&\leq 2 |\good| \exp\left(-\samples \frac{\eps^2/4}{C_1 + C_2 \eps/2}\right). \nonumber 
\end{align}
Similarly we can bound the probability of the bad components 
being higher than the threshold,
\begin{align}
\P(\max_{k\in\bad}&|\ip{Ay}{A\atom_k} |\geq \max_{k\in\bad}|\ip{y}{\atom_k} | +\frac{\eps}{2}) \notag\\
&\leq  \P(\bigcup_{k\in \bad}\{ |\ip{Ay}{A\atom_k}| \geq |\ip{y}{\atom_k}| + \frac{\eps}{2}\})\nonumber \\
&\leq \sum_{k\in\bad} \P( |\ip{Ay}{A\atom_k} - \ip{y}{\atom_k}| \geq \frac{\eps}{2})\notag\\
&\leq 2 |\bad| \exp\left(-\samples \frac{\eps^2/4}{C_1 + C_2 \eps/2}\right). \nonumber 
\end{align}
Combining these two estimates we see that the probability of success 
for thresholding is exceeding
\begin{align}
1 - 2 \natoms \exp\left(-\samples \frac{\eps^2/4}{C_1 + C_2 \eps/2}\right). \nonumber 
\end{align}
The lemma finally follows from requiring this probability to be higher 
than $1-e^{-t}$ and solving for $\samples$.
\end{Proof}

The result above may appear surprising because the number of measurements 
seems to be independent of the sparsity. 
The dependence, however, is quite well hidden in the margin $\eps$ and the 
normalization $\|y\|_2=1$. For clarification we will estimate $\eps$ given 
the coefficients and the coherence of the dictionary.

\begin{corollary} Let $\dico$ be an $\spacedim \times \natoms$ dictionary
with Babel function $\mu_1$ defined in (\ref{def:Babel}).
Assume a signal $y=\dico_\good x$ with $|\good|=\sparsity$ 
satisfies the sufficient recovery 
condition for thresholding, 
\begin{equation}\label{cond:thresh_recover}
\frac{|x_{\min}|}{\|x\|_\infty}>\mu_1(\sparsity)+\mu_1(\sparsity - 1),
\end{equation}
where $|x_{\min}| = \min_{i \in \good} |x_i|$. If $A$ is an $\samples \times \spacedim$ random matrix according to one of the probability models in Lemma 
\ref{lemma:ipconc} then with probability at least $1-e^{-t}$
thresholding can recover $x$ (and hence $y$) from $s = A y = A\dico x$ as long as
\begin{align}\label{cor:threshsamplesdico} 
\samples \geq& C_3 \sparsity (1+ \mu_1(\sparsity-1)) (\log(2\natoms)+t) \notag\\
&\cdot \left( \frac{|x_{\min}|}{\|x\|_\infty}-\mu_1(\sparsity)-\mu_1(\sparsity - 1)\right)^{-2}.
\end{align}
Here, $C_3$ is the constant from Theorem \ref{th:recth}.

In the special case that the dictionary is an ONB the signal 
always satisfies the recovery condition and the bound for the 
necessary number of samples reduces to
\begin{equation}\label{cor:threshsamplesONB} \samples > C_3  \sparsity \left(\frac{\|x\|_\infty}{|x_{\min}|}\right)^2 (\log(2\natoms)+t).
\end{equation}
\end{corollary} 
\begin{Proof}
%It is shown for instance in \cite{ICASSP} that the (\ref{cond:thresh_recover})
%is a sufficient recovery condition for thresholding (applied with $\dico$ on $y$).
%
The best possible value for $\eps$ in Theorem~\ref{th:recth} is quite obviously
\begin{align}
\eps &= \min_{i\in\good}|\ip{y/\|y\|_2}{\atom_i}|- \max_{k\in\bad}|\ip{y/\|y\|_2}{\atom_k} | \notag\\
&= \frac{1}{\|y\|_2}\big(| \min_{i\in\good} \sum_{j\in\good} x_j \ip{\atom_j}{\atom_i}| \notag\\
&\phantom{=\frac{1}{\|y\|_2}(| \min_{i\in\good} \sum_{j\in\good} x_j } \:- \max_{k \in \bad} |\sum_{j\in\good} x_j \ip{\atom_j}{\atom_k}|\big)\notag\\
&\geq  \frac{1}{\|y\|_2}\left( |x_{\min}| - \|x\|_\infty \mu_1(\sparsity-1) -\|x\|_\infty \mu_1(\sparsity)\right).
\notag
\end{align}
Therefore, we can bound the factor $C(\eps)$ in Theorem~\ref{th:recth} as 
\begin{align}
C&(\eps)\leq C_3 \eps^{-2}\notag\\
& \leq C_3  \frac{\|y\|_2^2}{\|x\|^2_\infty} \cdot 
\big( \frac{|x_{\min}|}{\|x\|_\infty}-\mu_1(\sparsity)-\mu_1(\sparsity - 1)\big)^{-2}.\notag
\end{align}
To get to the final estimate observe that by Lemma \ref{lem:coh}
\begin{align}
 \frac{\|y\|_2^2}{\|x\|^2_\infty} =\frac{\|\dico_\good x \|_2^2}{\|x\|^2_\infty}& \leq (1 + \mu_1(\sparsity-1)) \frac{\| x \|_2^2}{\|x\|^2_\infty} \notag\\
 & \leq  (1 + \mu_1(\sparsity-1)) \sparsity.\notag
\end{align}
The case of an ONB simply follows from $\mu_1(\sparsity)=0$.
\end{Proof}

The previous results tell us that as for BP 
we can choose the number $\samples$ 
of samples linear in the sparsity $\sparsity$. However, 
for thresholding successful recovery additionally depends 
on the ratio of the largest
to the smallest coefficient. Also, in contrast to BP the result 
is no longer
uniform, meaning that the stated success probability is only valid 
for the given 
signal $x$. It does not imply that a single matrix $A$ can ensure
recovery for all sparse signals. Indeed, in the case of a 
Gaussian matrix $A$ and
an orthonormal basis $\dico$ it is known that once $A$ is randomly chosen then 
with high probability there exists
a sparse signal $x$ (depending on $A$) such that thresholding fails on $x$ 
unless the number of samples
$\samples$ is quadratic in the sparsity $\sparsity$, 
see e.g.~\cite[Section 7]{do06}. This fact seems to generalise
to redundant $\dico$.
%There seems to be no reason why this situation should improve when
%we generalize from orthonormal bases to redundant dictionaries $\dico$.

%one cannot guarantee recovery for all
%signals of sparsity $\sparsity$ with high probability unless

\begin{example}[Dirac-DCT] 
Assume again that our dictionary is 
the union of the Dirac and the Discrete Cosine Transform bases in $\R^d$ for $d=2^{2p+1}$. The 
coherence is again $\mu=2^{-p}$ and the number of atoms $K= 2^{2p+1}$. 
If we assume the sparsity $S \leq 2^{p-2}$ 
and balanced coefficients, \ie $|x_i|=1$, we get the following crude 
estimate for the number of necessary samples
\begin{equation}\nonumber
\samples\geq 6C_3\,\sparsity (\log(2)(2p+2) + t).
\end{equation}
If we just allow the use of one of the two ONBs to build the 
signal, the number of necessary samples reduces to
\begin{equation}\nonumber
\samples \geq C_3\,\sparsity  (\log(2)(2p + 1)+t).
\end{equation}
\end{example}

Again we see that whenever the sparsity $\sparsity \lesssim \sqrt{d}$ the 
results for ONBs and general dictionaries are comparable. 
At this point it would be nice to have a similar result for OMP. 
This task seems rather difficult due to stochastic dependency issues and 
so, unfortunately, we have not been able to do this analysis yet. 
%Thus, the only thing remaining to be done is to see how things work in 'reality'. 

%%% Local Variables:
%%% TeX-master: "CompSensing.tex"
%%% End:
 
\section{Numerical Simulations}

For our numerical simulations we used the same dictionary as for the examples, \ie
the combination of the Dirac and the Discrete Cosine Transform bases 
in $\R^{\spacedim}$, $\spacedim=256$, with coherence 
$\mu=\sqrt{1/128}\approx 0.0884$.

\begin{figure}[htb] 
\centering

%%%%  for one paragraph
\begin{tabular}{cc}
 \includegraphics[height=2.2cm,width=5.5cm]{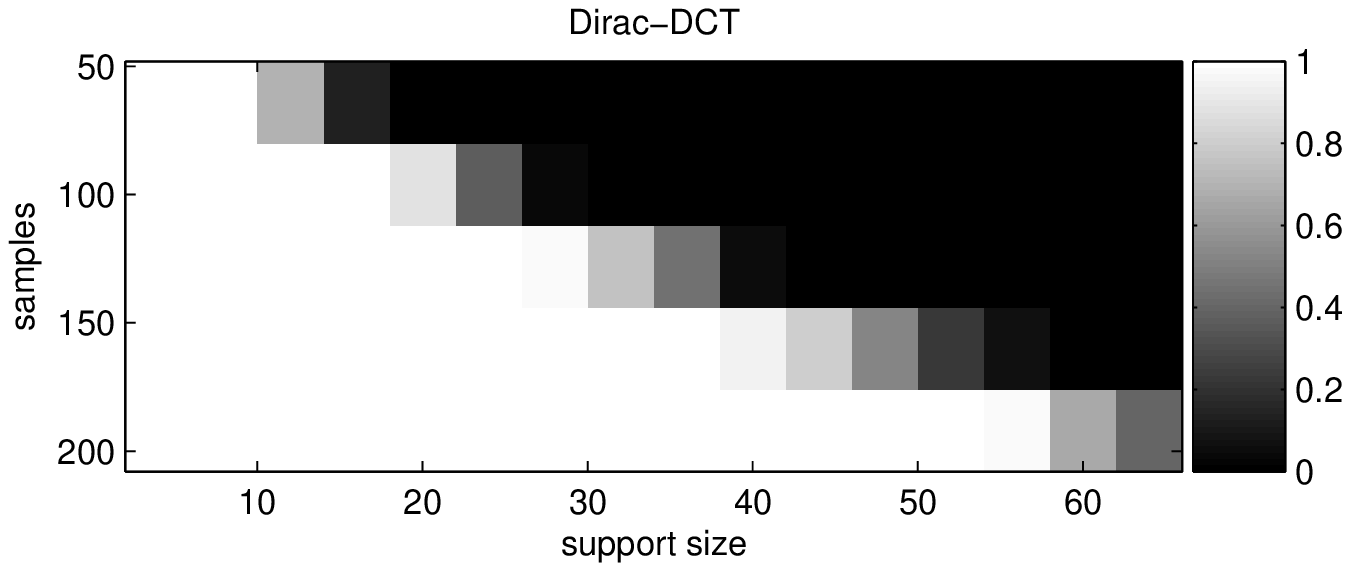} &
  \includegraphics[height=2.2cm,width=5.5cm]{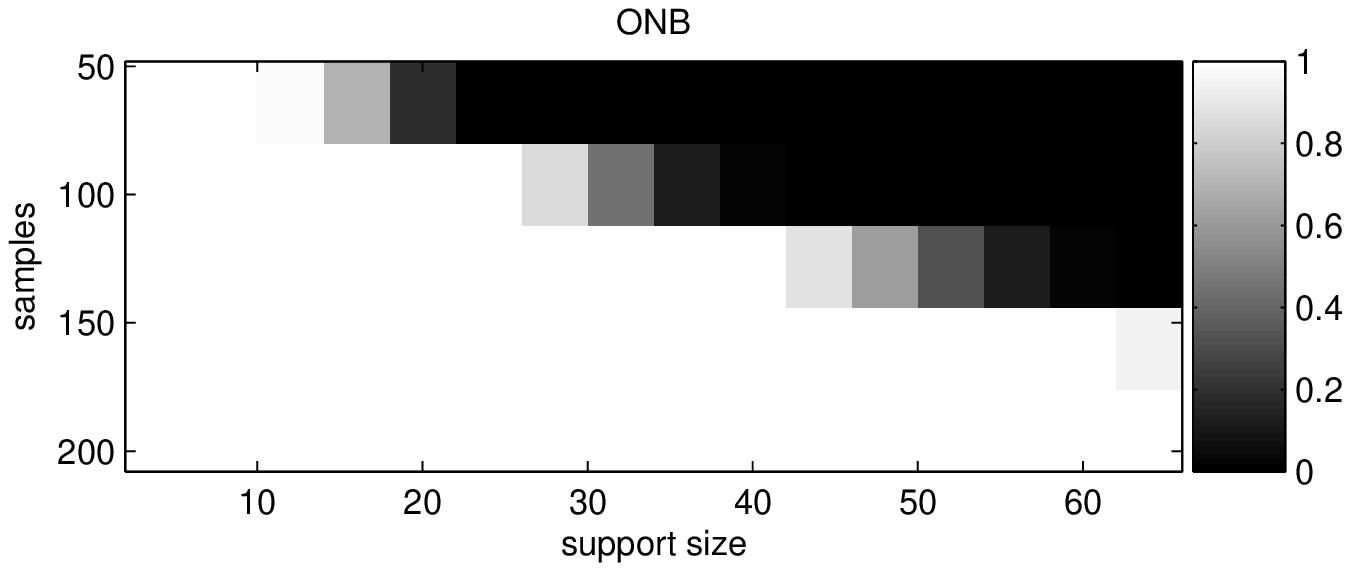}\\
\end{tabular}

%%%% for two paragraphs
%\begin{tabular}{c}
% \includegraphics[height=2.2cm,width=5.5cm]{figures/bp_diracdct.eps} \\
%  \includegraphics[height=2.2cm,width=5.5cm]{figures/bp_dirac.eps}\\
%\end{tabular}

\caption{Recovery Rates for BP as a Function of the Support and Sample Sizes \label{figure:bp}}
\end{figure}

\begin{figure}[htb] 
\centering

%%%%  for one paragraph
\begin{tabular}{cc}
 \includegraphics[height=2.5cm,width=5.5cm]{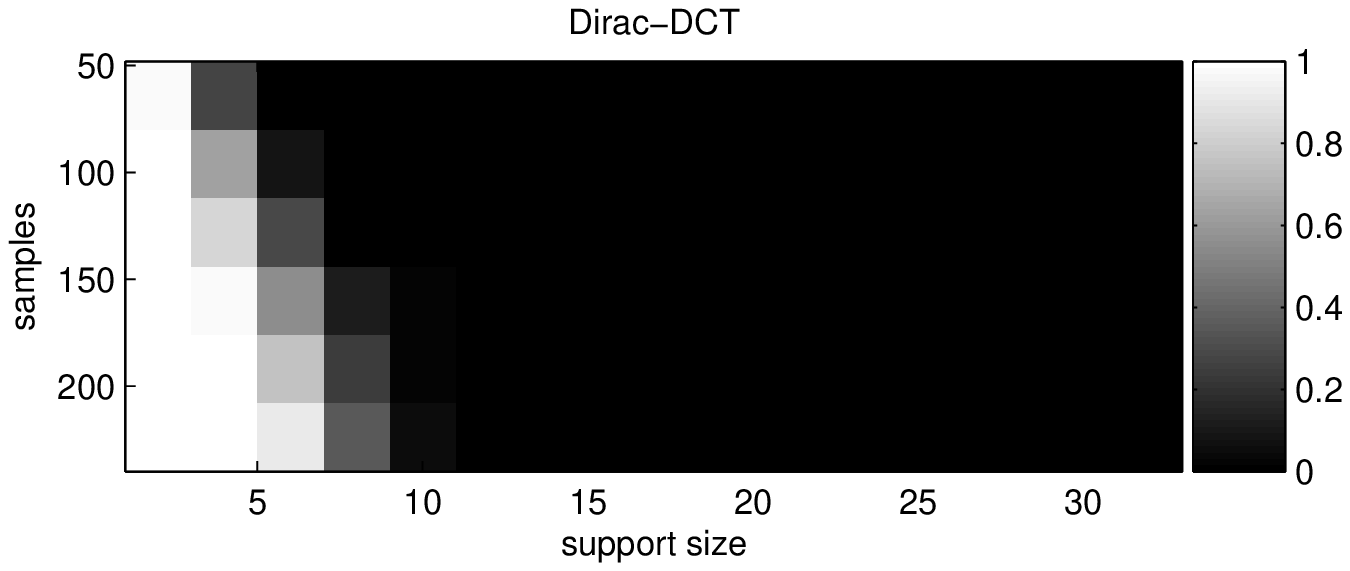} &
  \includegraphics[height=2.5cm,width=5.5cm]{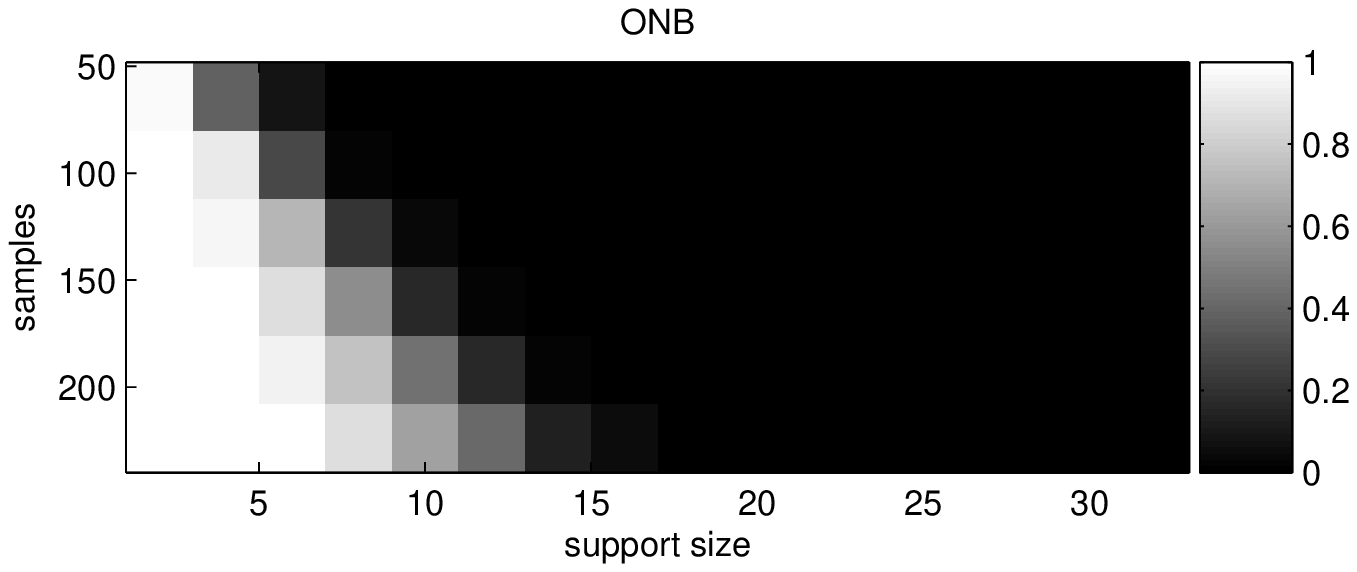}\\
\end{tabular}

%%%% for two paragraphs
%\begin{tabular}{c}
% \includegraphics[height=2.5cm,width=5.5cm]{figures/th_diracdct.eps} \\
%  \includegraphics[height=2.5cm,width=5.5cm]{figures/th_dirac.eps}\\
%\end{tabular}

      \caption{Recovery Rates for Thresholding  as a Function of the Support and Sample Sizes \label{figure:thresh}}
\end{figure}

\begin{figure}[htb] 
\centering

%%%%  for one paragraph
\begin{tabular}{cc}
 \includegraphics[height=2.5cm,width=5.5cm]{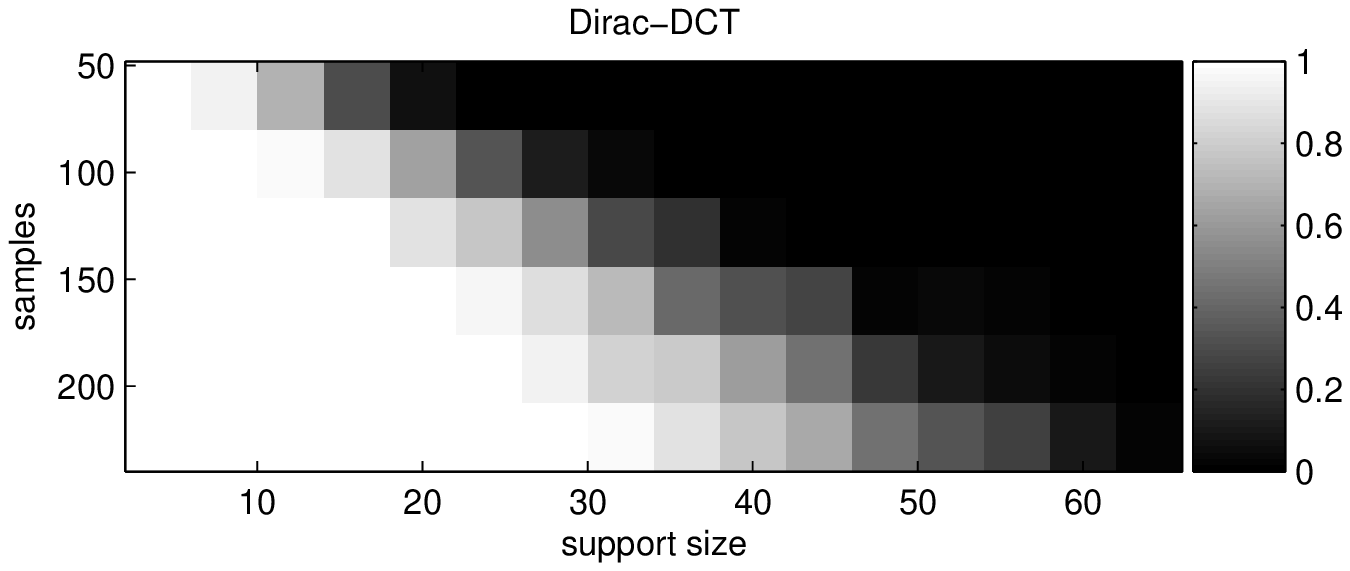} &
  \includegraphics[height=2.5cm,width=5.5cm]{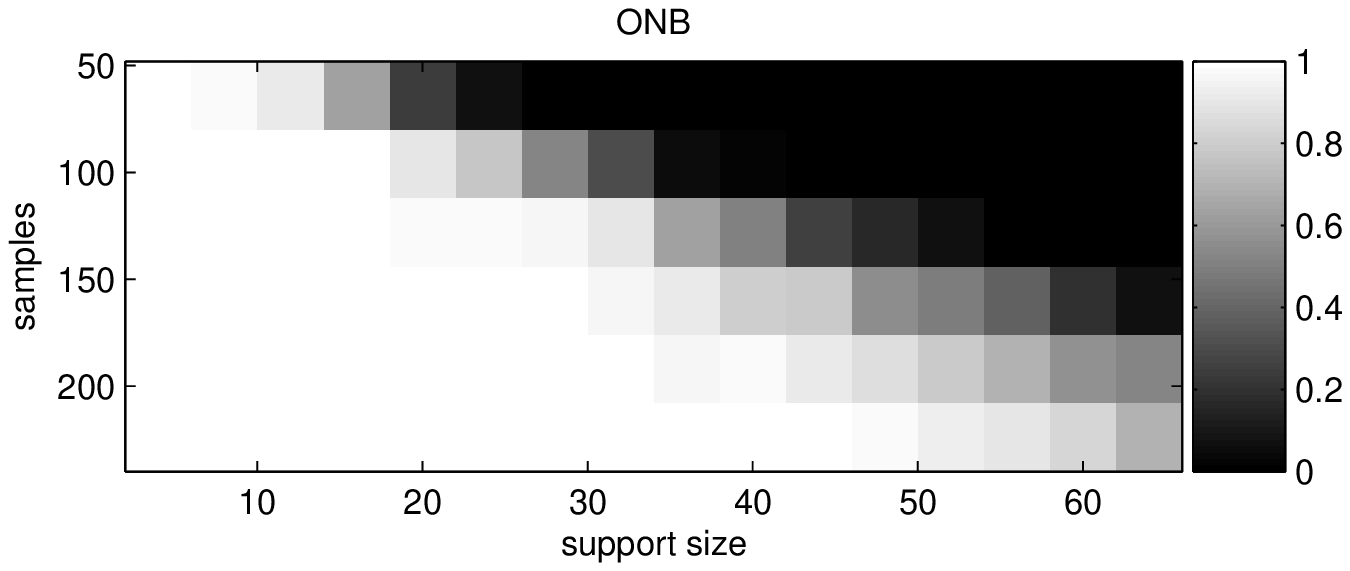}\\
\end{tabular}

%%%% for two paragraphs
%\begin{tabular}{c}
% \includegraphics[height=2.5cm,width=5.5cm]{figures/mp_diracdct.eps}\\
%  \includegraphics[height=2.5cm,width=5.5cm]{figures/mp_dirac.eps}\\
%\end{tabular}

      \caption{Recovery Rates for OMP  as a Function of the Support and Sample Sizes  \label{figure:omp}}
\end{figure}

We drew six measurement matrices of size $n\times \spacedim$, 
with $n$ varying between 64 and 224 in steps of 32, by choosing each 
entry as independent realisation of a centered Gaussian random variable 
with variance $\sigma^2=n^{-1}$. 
Then for every sparsity level $\sparsity$, varying between 4 and 64 
in steps of 4, respectively between 2 and 32 in steps of 2 for thresholding, 
we constructed 100 signals. The support $\good$ was chosen uniformly at 
random among all $\natoms \choose \sparsity$ possible supports of the 
given sparsity $\sparsity$. For BP and OMP
the coefficients $(x_i)_{i\in\good}$ of the 
corresponding entries were drawn from a normalised standard Gaussian 
distribution while for thresholding we chose them of absolute value one 
with random signs. Then for each of the algorithms we counted 
how often the correct support could be recovered. For comparison the same 
setup was repeated replacing the dictionary with the canonical (Dirac) basis. The 
results are displayed in Figures~\ref{figure:bp},~\ref{figure:thresh} 
and~\ref{figure:omp}.

%%%%%%%%%%%%%%%%%%%%%%%%%
%%%%%%%%%%%%%%%%%%%%%%%%%%%

%%%%%%%%%%%%%%%%%%%%%%%%%%%%
%%%%%%%%%%%%%%%%%%%%%%%%%%%%%%%

%%%%%%%%%%%%%%%%%%%%%%%%%%%%%
%%%%%%%%%%%%%%%%%%%%%%%%%%%%
As predicted by the theorems the necessary number of measurements is higher if the 
sparsity inducing dictionary is not an ONB. If we compare the 
three recovery schemes 
we see that thresholding gives the weakest results as expected. 
However, the improvement in performance of BP over OMP is not that significant. This is 
especially interesting considering that in practice BP is a lot more computationally 
intensive than OMP.

%%% Local Variables:
%%% TeX-master: "CompSensing.tex"
%%% End:

\section{Conclusions \& Future Work}

We have shown that compressed sensing can also be applied to signals that 
are sparse in a redundant dictionary. The spirit is that whenever the support 
can be reconstructed from the signal itself it can also be reconstructed 
from a small number of random samples with high probability. 
We have shown that this kind of 
stability is valid for reconstruction by Basis Pursuit as well as for 
the simple thresholding algorithm.
Thresholding has the advantage of being much faster and easier to implement than BP.
However, it has the slight drawback that the number of required samples depends 
on the ratio of the largest to the smallest coefficient, and
recovery is only guaranteed with high probability for a given signal and not
uniformly for all signals in contrast to BP. 
Furthermore, there is numerical evidence that Orthogonal Matching Pursuit
also works well. In particular, it is still faster than BP and
the required number of samples does not seem to depend on the ratio of the largest
to the smallest coefficient.

For the future there remains plenty of work to do. First of all we 
would like to have a recovery theorem for OMP comparable to Theorem 
\ref{th:recth}. However, since in the course of iterating 
the updated residuals become stochastically dependent on the random matrix $A$ this
task does not seem to be straightforward. In particular, the technique developed in 
\cite{gitr05} cannot be applied directly.
Then we would like to investigate for which dictionaries it is possible
to replace the random Gaussian/Bernoulli matrix by a random Fourier matrix, 
see also \cite{ra05-7}.
This would have the advantage that the Fast Fourier Transform can be used in
the algorithms in order to speed up the reconstruction. 
Finally, it would be interesting to relax the incoherence assumption on the 
dictionary. 

%%% Local Variables:
%%% TeX-master: "CompSensing.tex"
%%% End:

\appendix
\section{Proof of Lemma \ref{lemma:ipconc}} \label{app:proof}

%In this section we investigate recovery from random
%measurements by OMP of signals
%that are sparse with respect to redundant dictionaries.
%We proceed by generalizing results by Gilbert and
%Tropp in \cite{gitr05}. 

Our proof uses the following inequality due to Bennett 
(also refered to as Bernstein's inequality)
\cite[eq.~(7)]{be62}, see also \cite[Lemma~2.2.11]{vawe96}.
\begin{Theorem}\label{thm:Bennett} Let $X_1,\hdots,X_\samples$ be independent
random variables with zero mean such that 
\begin{equation}\label{moment_bound}
\E|X_i|^q \leq q! M^{q-2} v_i / 2
\end{equation}
for every $m \geq 2$ and some constants $M$ and $v_i$, $i = 1,\hdots,\samples$.
Then for $x > 0$
\[
\P\left(|\sum_{i=1}^\samples X_i| \geq x\right) \leq 2e^{-\frac{1}{2}\frac{x^2}{v + M x}}
\]
%and
%\[
%\P\left(\sum_{i=1}^\samples X_i\right \leq -x\right) \leq e^{-\frac{1}{2}\frac{x^2}{v + M x}}
%\]
with $v = \sum_{i=1}^\samples v_i$.
\end{Theorem}

Now let us prove Lemma \ref{lemma:ipconc}.
Observe that
\[
\langle A x, Ay \rangle = 
\frac{1}{\samples} \sum_{\ell=1}^\samples \sum_{k=1}^\spacedim \sum_{j=1}^\spacedim
g_{\ell k} g_{\ell j} x_k y_j
\]
where $g_{\ell k}$, $\ell=1,\hdots,\samples, k=1,\hdots,\spacedim$ are 
independent
standard Gaussians. We define the random variable
\[
Y:= \sum_{k,j=1}^\spacedim g_k g_j x_k y_j
\]
where again the $g_k$, $k=1,\hdots,\spacedim$ are independent standard 
Gaussians.
Then we can write
\[
\langle Ax, Ay \rangle= \frac{1}{\samples} \sum_{\ell=1}^\samples Y_\ell
\]
where the $Y_\ell$ are independent copies of $Y$.

Let us investigate $Y$. The expectation of $Y$ is easily calculated as
\[
\E Y =\sum_{k=1}^\spacedim x_k y_k = \langle x,y\rangle.
\]
Hence, also $\E \left[ \langle A x, A y\rangle \right]= \langle x,y\rangle$.
Now let
\[
Z:= Y - \E Y = \sum_{k \neq j} g_j g_k x_j x_k + \sum_{k} (g_k^2 - 1) x_k y_k.
\]
The random variable $Z$ is known as Gaussian chaos of order $2$.

Thus, we have to show the moment bound \eqref{moment_bound} for the random
variable $Z$. Note that $\E Z = 0$. 
A general bound for Gaussian chaos (see \cite[p.~65]{leta91})
gives
\begin{equation}\label{chaos:estimate}
\E |Z|^q \leq (q-1)^q \left( \E |Z|^2\right)^{q/2}
\end{equation}
for all $q \geq 2$. 
Using Stirling's formula, $q! = \sqrt{2\pi q}\, q^q e^{-q} e^{R_q}$, $\frac{1}{12q+1} \leq R_q \leq \frac{1}{12q}$, we further obtain, for all $q \geq 3$:
%\begin{align}
%\E |Z|^q &= q! \frac{(q-1)^q}{e^{R_q} \sqrt{2\pi q}\, e^{-q}q^{q}} 
%\left(\E |Z|^2\right)^{q/2}\notag\\
%&= \left(1-\frac{1}{q}\right)^q \frac{e^2 q! }{e^{R_q} \sqrt{2\pi q}} \left(e^2 \E |Z|^2\right)^{\frac{q-2}{2}} \E |Z|^2\notag\\
%&\leq \frac{e}{e^{R_q} \sqrt{2\pi q}} q! \left(e^2 \E |Z|^2\right)^{\frac{q-2}{2}} 
%\E |Z|^2\notag\\
%&\leq q! \left( e (\E |Z|^2)^{1/2} \right)^{q-2} \frac{e}{\sqrt{6\pi}} \E |Z|^2. 
%\notag 
%\end{align}
\begin{align}
\E |&Z|^q = q! \frac{(q-1)^q}{e^{R_q} \sqrt{2\pi q}\, e^{-q}q^{q}} 
\left(\E |Z|^2\right)^{q/2}\notag\\
&= \left(1-\frac{1}{q}\right)^q \frac{e^2 q! }{e^{R_q} \sqrt{2\pi q}} \left(e^2 \E |Z|^2\right)^{(q-2)/2} \E |Z|^2\notag\\
&\leq \frac{e}{e^{R_q} \sqrt{2\pi q}} q! \left(e^2 \E |Z|^2\right)^{(q-2)/2} 
\E |Z|^2\notag\\
&\leq q! \left( e (\E |Z|^2)^{1/2} \right)^{q-2} \frac{e}{\sqrt{6\pi}} \E |Z|^2. 
\notag 
\end{align}

%A finer inspection shows that $C_1 \sqrt{q} \geq 1$ for all $q \geq 2$.
Hence, the moment bound \eqref{moment_bound} holds for all $q \geq 3$ with
\[ 
M = e \left(\E|Z|^2\right)^{1/2}, \qquad v = \frac{2e}{\sqrt{6\pi}} \E |Z|^2,
\]
and by direct inspection it then also holds for $q=2$.
So let us determine $\E |Z|^2$. Using independence of the $g_k$ we obtain
\begin{align}
\E |&Z|^2 =\E \left[ \sum_{j\neq k} \sum_{j' \neq k'} g_j g_k g_{j'} g_{k'}
x_j y_k x_{j'} y_{k'} \notag \right.\\
& \phantom{= \E [[} \left. +2 \sum_{j \neq k} \sum_{k'} g_j g_k (g_{k'}^2-1) x_j y_k x_{k'} y_{k'}\notag\right.\\
& \phantom{= \E [[} \left. 
+ \sum_{k} \sum_{k'} (g_k^2 - 1)(g_{k'}^2-1) x_k y_k x_{k'} y_{k'}\right]\notag\\
&= \sum_{k \neq j} \E [g_j^2] \E[g_k^2] x_jy_jx_k y_k \notag \\
& \hspace{1.5cm}+ \sum_{k \neq j} \E [g_j^2] \E[g_k^2] x_j^2 y_k^2 \notag\\
& \hspace{2cm} + \sum_{k} \E[(g_k^2-1)^2] x_k^2 y_k^2\label{eq:bernoulli} \\
&= \sum_{k \neq j}x_jy_jx_k y_k + \sum_{k \neq j} x_j^2 y_k^2 + 2 \sum_{k} x_k^2 y_k^2 \notag \\
&= \sum_{j,k}x_jy_jx_k y_k + \sum_{j,k} x_j^2 y_k^2 \notag \\
&=  \langle x,y\rangle^2+\|x\|_2^2 \|y\|_2^2 
\,\leq \, 2 \label{EZ2:estimate}
\end{align}
since by assumption $\|x\|_2, \|y\|_2 \leq 1$. Denoting by $Z_\ell$, $\ell = 1,\hdots,\samples$ independent copies of $Z$, Theorem \ref{thm:Bennett} yields
\begin{align}
\P\big( |\langle Ax,Ay \rangle& - \langle x,y\rangle| \geq t\big) \notag \\
&= \P\left(|\sum_{\ell=1}^\samples Z_\ell| \geq nt\right) \notag\\
&\leq 2e^{-\frac{1}{2} \frac{n^2 t^2}{nv + nMt}} =  2e^{-n \frac{t^2}{C_1 + C_2 t}} \notag,
\end{align}
with $C_1 =  \frac{2e}{\sqrt{6\pi}} E|Z|^2 \leq \frac{4e}{\sqrt{6\pi}}
\approx 2.5044$ and $C_2 =  e \sqrt{2} \approx 3.8442$. 
%Applying
%Theorem \ref{thm:Bennett} to $-Z$ yields the second inequality (\ref{sec:ineq})
%in Lemma \ref{lemma:ipconc}.

For the case of Bernoulli random matrices the proof is completely analogue. We just have to replace the standard Gaussians $g_k$ by $\pm 1$ Bernoulli variables.
In particular, the estimate (\ref{chaos:estimate}) for the chaos variable $Z$
is still valid, see \cite[p.~105]{leta91}. Furthermore, for Bernoulli 
variables $g_k$ we clearly have $g_k^2 = 1$.
Hence, going through the estimate above we see that in \eqref{eq:bernoulli} the last
term is actually zero, so the final bound in \eqref{EZ2:estimate} is still valid.

%In particular, we remove their
%assumption that the columns of the measurement matrix
%be stochastically independent. This is even interesting for
%the case of sparsity with respect to the canonical bases.

%%% Local Variables:
%%% TeX-master: "CompSensing.tex"
%%% End:
 
\bibliographystyle{abbrv}
\bibliography{CompSensingBib}

\begin{thebibliography}{10}

\bibitem{ac01}
D.~Achlioptas.
\newblock Database-friendly random projections.
\newblock In {\em Proc. 20th Annual ACM SIGACT-SIGMOD-SIGART Symp. on
  Principles of Database Systems}, pages 274--281, 2001.

\bibitem{ahelbr06}
M.~Aharon, M.~Elad, and A.~Bruckstein.
\newblock {K}-{S}{V}{D}: An algorithm for designing of overcomplete
  dictionaries for sparse representation.
\newblock {\em IEEE Trans. on Signal Processing.}, 54(11):4311--4322, November
  2006.

\bibitem{badadewa06}
R.~{B}araniuk, M.~{D}avenport, R.~{D}e{V}ore, and M.~{W}akin.
\newblock {A} simple proof of the restricted isometry property for random
  matrices.
\newblock {\em {C}onstr. {A}pprox.}, to appear.

\bibitem{be62}
G.~{B}ennett.
\newblock {P}robability inequalities for the sum of independent random
  variables.
\newblock {\em {J}. {A}m. {S}tat. {A}ssoc.}, 57:33--45, 1962.

\bibitem{carota06-1}
E.~{C}and{\`e}s, J.~{R}omberg, and T.~{T}ao.
\newblock {S}table signal recovery from incomplete and inaccurate measurements.
\newblock {\em {C}omm. {P}ure {A}ppl. {M}ath.}, 59(8):1207--1223, 2006.

\bibitem{cataXX}
E.~{C}and{\`e}s and T.~{T}ao.
\newblock {N}ear optimal signal recovery from random projections: universal
  encoding strategies?
\newblock {\em {I}{E}{E}{E} {T}rans. {I}nf. {T}heory}, 52(12):5406--5425, 2006.

\bibitem{chdosa99}
S.~{C}hen, D.~{D}onoho, and M.~{S}aunders.
\newblock {A}tomic decomposition by {B}asis {P}ursuit.
\newblock {\em {S}{I}{A}{M} {J}. {S}ci. {C}omput.}, 20(1):33--61, 1999.

\bibitem{avdama97}
G.~{D}avis, S.~{M}allat, and M.~{A}vellaneda.
\newblock {A}daptive greedy approximations.
\newblock {\em {C}onstr. {A}pprox.}, 13(1):57--98, 1997.

\bibitem{do04}
D.~{D}onoho.
\newblock {C}ompressed {S}ensing.
\newblock {\em {I}{E}{E}{E} {T}rans. {I}nf. {T}heory}, 52(4):1289--1306, 2006.

\bibitem{doelte06}
D.~{D}onoho, M.~{E}lad, and V.~{T}emlyakov.
\newblock {S}table recovery of sparse overcomplete representations in the
  presence of noise.
\newblock {\em {I}{E}{E}{E} {T}rans. {I}nf. {T}heory}, 52(1):6--18, 2006.

\bibitem{dota06}
D.~{D}onoho and J.~{T}anner.
\newblock {C}ounting faces of randomly-projected polytopes when the projection
  radically lowers dimension.
\newblock {\em Preprint arXiv:math.MG/0607364}, 2006.

\bibitem{do06}
D.~L. {D}onoho.
\newblock {F}or most large underdetermined systems of linear equations the
  minimal $\ell^1$-norm solution is also the sparsest solution.
\newblock {\em {C}omm. {P}ure {A}ppl. {M}ath.}, 59(6):797--829, 2006.

\bibitem{dr06}
I.~Drori.
\newblock Fast $\ell_1$ minimization by iterative thresholding for
  multidimensional {N}{M}{R} spectroscopy.
\newblock {\em Preprint}, 2006.

\bibitem{gitr05}
A.~C. {G}ilbert and J.~A. {T}ropp.
\newblock {S}ignal recovery from random measurements via orthogonal matching
  pursuit.
\newblock {\em {I}{E}{E}{E} {T}rans. {I}nform. {T}heory}, to appear.

\bibitem{grmarascva06}
R.~{G}ribonval, B.~{M}ailhe, H.~{R}auhut, K.~{S}chnass, and P.~{V}andergheynst.
\newblock {A}verage case analysis of multichannel thresholding.
\newblock In {\em Proc. IEEE ICASSP07, Honolulu}, 2007.

\bibitem{leta91}
M.~{L}edoux and M.~{T}alagrand.
\newblock {\em {P}robability in {B}anach spaces. {I}soperimetry and processes.}
\newblock {S}pringer-{V}erlag, {B}erlin, {H}eidelberg, {N}ew{Y}ork, 1991.

\bibitem{mazh93}
S.~G. {M}allat and Z.~{Z}hang.
\newblock {M}atching pursuits with time-frequency dictionaries.
\newblock {\em {I}{E}{E}{E} {T}rans. {S}ignal {P}rocess.}, 41(12):3397--3415,
  1993.

\bibitem{mepato06}
S.~{M}endelson, A.~{P}ajor, and N.~{T}omczak{-}{J}aegermann.
\newblock {U}niform uncertainty principle for {B}ernoulli and subgaussian
  ensembles.
\newblock {\em Preprint}, 2006.

\bibitem{na95}
B.~{N}atarajan.
\newblock {S}parse approximate solutions to linear systems.
\newblock {\em {S}{I}{A}{M} {J}. {C}omput.}, 24:227--234, 1995.

\bibitem{ra05-7}
H.~{R}auhut.
\newblock {R}andom sampling of sparse trigonometric polynomials.
\newblock {\em {A}ppl. {C}omput. {H}arm. {A}nal.}, 22(1):16--42, 2007.

\bibitem{ru06-1}
M.~{R}udelson and R.~{V}ershynin.
\newblock {S}parse reconstruction by convex relaxation: {F}ourier and
  {G}aussian measurements.
\newblock In {\em {P}roc. {C}{I}{S}{S} 2006 (40th {A}nnual {C}onference on
  {I}nformation {S}ciences and {S}ystems)}, 2006.

\bibitem{tr04}
J.~{T}ropp.
\newblock {G}reed is good: {A}lgorithmic results for sparse approximation.
\newblock {\em {I}{E}{E}{E} {T}rans. {I}nf. {T}heory}, 50(10):2231--2242, 2004.

\bibitem{vawe96}
A.~{V}an~der {V}aart and J.~{W}ellner.
\newblock {\em {W}eak convergence and empirical processes}.
\newblock {S}pringer-{V}erlag, 1996.

\end{thebibliography}

\end{document}